\documentclass[12pt]{article}
\usepackage{amsthm,amsfonts,amsmath}

\input epsf.tex
\newdimen\epsfxsize
\newdimen\epsfysize
\def\qed{\vrule height5pt width3pt depth.5pt}

\theoremstyle{plain}
\newtheorem{thm}{Theorem}[section]

\newtheorem{prop}[thm]{Proposition}
\newtheorem{conj}[thm]{Conjecture}

\newtheorem{rem}{Remark}[section]

\begin{document}

\title{Non-Trivial Realizations of Virtual Link Diagrams}         
\author{H. A. Dye \\
MADN-MATH \\
United States Military Academy \\
646 Swift Road \\
West Point, NY 10996 \\
hdye@ttocs.org }        
\date{February 23, 2005}          
\maketitle

\begin{abstract}
A realization of a virtual link diagram is obtained by choosing over/under markings for each virtual crossing. Any realization can also be obtained from some representation of the virtual link. (A representation of a virtual link is a link diagram on an oriented 2-dimensional surface.) We prove that if a minimal genus representation meets certain criteria then there is a minimal genus representation resulting in a knotted realization. 
\end{abstract}

\noindent {\bf Acknowledgement.}The views expressed herein are those of the author and do not purport to reflect the position of the United States Military Academy, the Department of the Army, or the Department of Defense.
 (Copyright
2005.)

\section{Introduction}

We can construct a realization of any virtual link diagram by choosing  over/under markings for each virtual crossing. This process produces a set of classical link diagrams. A virtual link diagram with $ n $ virtual crossings has $ 2^n $ possible realizations. 

Realizations were originally investigated by Sam Nelson \cite{nelson}. In this paper, Nelson showed that a sequence of virtual Reidemeister and Reidemeister moves performed on a virtual link diagram could not always be replicated by a sequence of classical Reidemeister moves performed on a realization.

A realization of a virtual link diagram in the equivalence class of the virtual link $ \hat{L} $ can be obtained from a representation of $ \hat{L} $.
However, the realizations of a virtual link form a set of inequivalent classical link diagrams. This provokes the question: What is the relationship between two realizations of a virtual link?

In \emph{Unsolved Problems in Virtual Knot Theory} \cite{vkprobs}, the authors ask: Does every virtual link  have a non-trivial realization obtained from an unknotted minimal genus representation?
Relaxing the restrictions that the representation have an unknotted and minimal genus surface,
we can easily construct knotted realizations of any virtual link.  
A non-minimal representation with an unknotted surface can be reduced to a minimal genus surface via a sequence of handle cancellations, Reidemeister moves in the surface, and homeomorphisms of the surface \cite{kup}.
Unfortunately, there is no guarantee that the knotting or linking in the realization survives a sequence of handle cancellations and homeomorphisms of the surface.  
For example, homemorphisms of the surface can  unknot or unlink the corresponding realizations. Knotting of the representation's surface can be induced by handle cancellation or addition.

An alternative approach to the problem is to start with a minimal genus representation $ (F,L) $ such that $ F $ is an unknotted abstract surface.   The classical diagram $ L $ in $ \mathbb{R}^3 $ (instead of on the surface $ F $) forms a realization of a virtual link diagram, $ \hat{L} $. (We recover $ \hat{L} $ by projecting $ L $ onto the plane while preserving the over/under markings at each classical crossing and marking any new crossings as virtual.) If the realization is trivial, we want a sequence of moves that produce an equivalent representation resulting in a non-trivial realization. These moves can not increase the genus of the surface or knot the surface.

The standard diagram of Kishino's knot \cite{kishpoly} illustrates the difficulties of finding a non-trivial realization that is obtained from a minimal genus representation. The common minimal genus representations of this knot result in trivial realizations. These realizations correspond to realizations of the standard diagram. 

In this paper, we demonstrate that if a minimal genus representation meets certain criteria then we may obtain a non-trivial realization. We begin with a minimal genus representation with an unknotted abstract surface that corresponds to a trivial realization ( $L $ viewed in $ \mathbb{R}^3 $). If this representation meets the criteria, we apply a sequence of Reidemeister moves (to the link) in the surface and perform Dehn twists on the abstract surface. This process produces a minimal genus representation with an unknotted surface that corresponds to a non-trivial realization.

\section{Virtual Links}

\begin{figure}[htb] \epsfysize = 1 in
\centerline{\epsffile{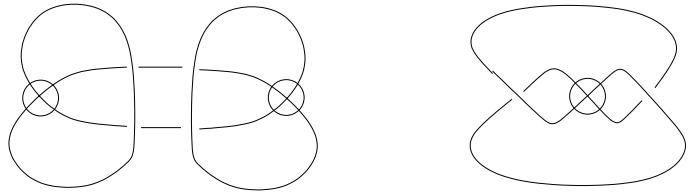}}
\caption{Virtual Link Diagrams}
\label{fig:vlinks}
\end{figure}

A virtual link diagram is a decorated immersion of $ n $ copies of $ S^1 $ in the plane. The diagram contains two types of crossings: classical crossings and virtual crossings. Classical crossings are indicated by under/over markings. Virtual crossings are indicated by a solid encircled X. Note that the classical link diagrams are a subset of the virtual link diagrams.
Two virtual link diagrams are shown in figure \ref{fig:vlinks}.

\begin{figure}[htb] \epsfysize = .75 in
\centerline{\epsffile{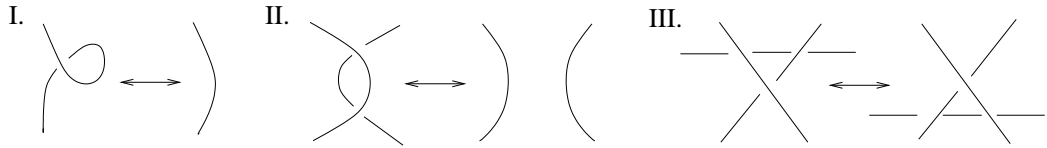}}
\caption{Reidemeister Moves}
\label{fig:rmoves}
\end{figure}

The \emph{Reidemeister moves} and the \emph{virtual Reidemeister moves} establish equivalence classes of diagrams. The Reidemeister moves involve only classical crossings and are shown in figure \ref{fig:rmoves}.

\begin{figure}[htb] \epsfysize = 1 in
\centerline{\epsffile{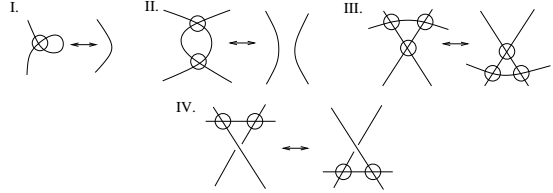}}
\caption{Virtual Reidemeister Moves}
\label{fig:vrmoves}
\end{figure}

The \emph{framed Reidemeister moves} refer to the Reidemeister II and III moves. 
The virtual Reidemeister moves are illustrated in figure \ref{fig:vrmoves}. Only the fourth virtual Reidemeister move involves classical and virtual crossings.

Two virtual link diagrams are \emph{equivalent} if one diagram may be transformed in the other via a sequence of Reidemeister and virtual Reidemeister moves. 
A \emph{virtual link} is an equivalence class of equivalent virtual link diagrams.

\begin{figure}[htb] \epsfysize = 1 in
\centerline{\epsffile{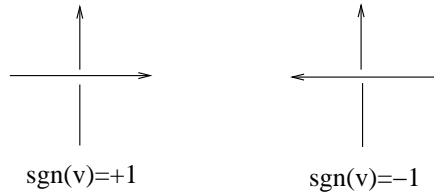}}
\caption{Crossing Sign}
\label{fig:xsign}
\end{figure}

We recall the definition of \emph{crossing sign}. We assign a value of $ \pm 1 $ to each classical crossing  as shown in figure \ref{fig:xsign}. 
The crossing sign of a classical crossing, $ v $, is denoted $ sgn(v) $.

 The \emph{writhe} of a virtual link diagram $ L $ is the sum of all crossing signs in the diagram. We denote the writhe of $ L $ as $ w(L)$:
\begin{equation}
w(L) = \underset{v \in L }{ \sum } sgn(v)
\end{equation}
The writhe is invariant under the framed Reidemeister moves and the virtual Reidemeister moves. 

Let $ L $ be an $n $ component virtual link diagram, with components $ L_1, L_2, \dots L_n $.
We define the linking number of the components $L_i $ and $L_j $, denoted $ lk(L_i, L_j) $.
Now:
\begin{equation*}
lk(L_i,L_j) = \underset{ v \in L_i,L_j}{\sum} sgn(v)
\end{equation*}
Note that $ lk (L_i,L_j) = lk (L_j,L_i) $and that linking number is invariant under the classical and virtual Reidemeister moves. 
Two unlinked components have linking number zero. 
\begin{rem}
An alternative definition of linking number for virtual links given in \cite{gpv}. In this definition, $ lk (L_j, L_i) \neq lk (L_i, L_j) $.
\end{rem}

\begin{figure}[htb] \epsfysize = 1 in
\centerline{\epsffile{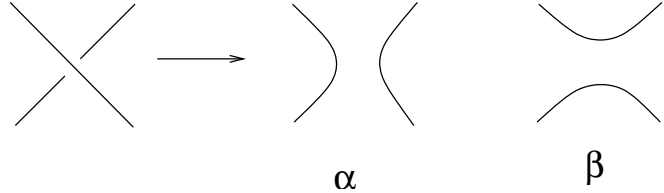}}
\caption{Smoothing Types}
\label{fig:smooth}
\end{figure}

We recall the Jones polynomial of a virtual link diagram. A \emph{smoothing} of a classical crossing removes a small neighborhood the diagram at the crossing. The crossed segments of the diagram are replaced with two non-intersecting segments. We  smooth a crossing horizontally (a type $ \alpha $ smoothing) or vertically (a type $ \beta $ smoothing) as shown in figure \ref{fig:smooth}.

In a \emph{state} of a virtual link diagram, each classical crossing is smoothed and implicitly labeled with its smoothing type.

We define the bracket polynomial of a virtual link diagram $ L $ (denoted as $ \langle L \rangle $). Let $ S $ represent the set of all possible states of $ L $ and let $ s $ denotes an element of $ S $. For a state $s$:    
\begin{gather*}
 |s| \text{ denotes the number of closed curves in the state } \\ \alpha_s \text{ represents the number of type $ \alpha $ smoothings,} \\  
\beta_s \text{ represents the number of type $ \beta $ smoothings.} 
\end{gather*}
Let $ d = -A^2 - A^{-2} $ then
\begin{equation*}
\langle L \rangle = \underset{s \in S}{ \sum } A^{\beta_s - \alpha_s} d^{ |s| -1 }
\end{equation*}
A state of a virtual link diagram consists of closed curves. (A state of a classical link diagram consists of simpled closed curves.) These states (possibly) contain virtual crossings.

The bracket polynomial is invariant under the framed Reidemeister and the virtual Reidemeister moves. 
The \emph{Jones polynomial} of a virtual link diagram $ L $ as
\begin{equation*} 
V(L) = \langle L \rangle A^{-3 w(L) }.
\end{equation*}
Recall that $ V(L) $ is invariant under the classical Reidemeister moves and the virtual Reidemeister moves \cite{kvirt}.

\begin{figure}[htb] \epsfysize = 0.75 in
\centerline{\epsffile{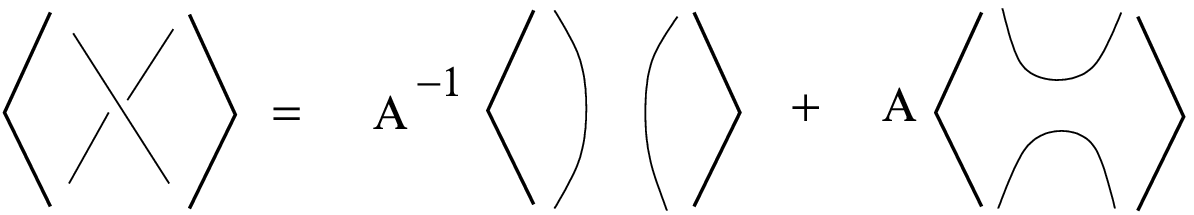}}
\caption{Skein Relation}
\label{fig:crossexpan}
\end{figure}

We may also compute the bracket polynomial by applying the skein relation, shown in figure \ref{fig:crossexpan}.

The skein relation is used to calculate equation \ref{jskein} given in Section \ref{proofs}.

\begin{figure}[htb] \epsfysize = 1.5 in
\centerline{\epsffile{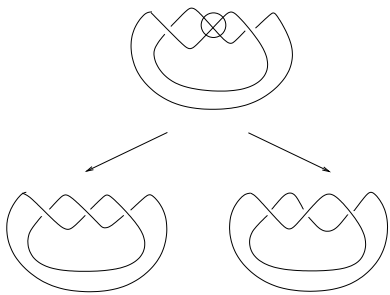}}
\caption{Virtual Knot and its Realizations}
\label{fig:vtrefreal}
\end{figure}

\begin{figure}[htb] \epsfysize = 1.5 in
\centerline{\epsffile{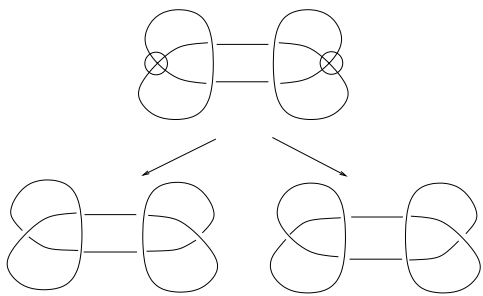}}
\caption{Kishino's Knot and Realizations}
\label{fig:kishreal}
\end{figure}

A \emph{realization} of a virtual link diagram with $ m $ virtual crossings is an assignment of over/under markings to each virtual crossing. For a virtual link diagram with $ m $ virtual crossings, there are $ 2^m $ realizations. A realization of a virtual link diagram with $n $ components is \emph{trivial} if the realization is equivalent to $ n $ unlinked copies of the unknot. We show some examples of virtual link diagrams and their realizations in figures \ref{fig:vtrefreal} and \ref{fig:kishreal}.

A virtual knot diagram with one virtual crossing and its two possible realizations are illustrated in figure \ref{fig:vtrefreal}. The realization on the right forms a trefoil which is not equivalent to the realization on the left(an unknot). A realization of a virtual link diagram with $ n $ components is \emph{trivial} if the realization is equivalent to $ n $ unknotted and unlinked copies of $ S^1 $.

The standard diagram of Kishino's knot has two virtual crossings and a set of four possible realizations.
Kishino's knot and two possible realizations are illustrated in figure \ref{fig:kishreal}. These realizations are equivalent to the unknot. 

Kishino's knot demonstrates that there exist virtual link diagrams with only trivial realizations.

Let $ \hat{L} $ be a virtual link diagram with $n $ classical crossings. We prove there exists a realization $ L $ of $ \hat{L} $ requiring $ n $ or fewer crossing changes to obtain the unlink.

\subsection{Unknotting Realizations}
We define the unknotting number of a fixed link diagram. 
The \emph{fixed unknotting number} of a fixed classical link diagram $ L $ is denoted as $ uf(L) $. This is the minimum number of crossing changes required to convert the fixed diagram $ L $ into the unlink. The fixed unknotting number can be used to define
the unknotting number \cite{adams} of a classical link $ L $. The unknotting number of $ L $ is denoted $ u(L) $ and $ u(L) = min \lbrace uf( \hat{L} ) | L \text{ is equivalent to } \hat{L} \rbrace $.

\begin{prop}
Let $ \hat{L} $ be a virtual link with $n $ classical crossings. Then there is a realization $ L $ such that $ uf(L) \leq n $ and all crossing changes occur at classical crossings.
\end{prop}
  
 \textbf{Proof:}
Let $ \hat{L} $ be a virtual link diagram with a realization $ L $. Suppose that $ uf(L) =g $ and that this realization has the minimal unknotting number of all realizations. If one of the crossing changes involves a virtual crossing then there exists a realization $ L_{r} $ such that $ uf(L_{r}) =g-1 $. This contradicts our assumption that $ L $ was minimal. As a result, $ uf(L) \leq n $ since the crossing changes occur only at classical crossings. \qed

\begin{figure}[htb] \epsfysize = 1 in
\centerline{\epsffile{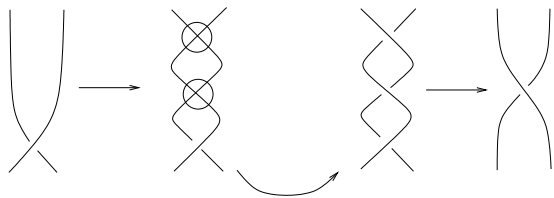}}
\caption{Obtaining an Unknotted Realization}
\label{fig:unknotreal}
\end{figure}

We may use this proposition to show:
\begin{prop}Let $ \hat{L} $ be a virtual link diagram with $ n $ classical crossings. Then there is an equivalent virtual link diagram related by a sequence of at most $ n $ virtual Reidemeister II moves that has a trivial realization.
\end{prop}
 \textbf{Proof:}
Let $ \hat{L} $ be a virtual link diagram with $ m $ virtual crossings and $ n $ classical crossings. By the above proposition, there exists a realization $ L $ with $ uf(L ) \leq n $. Identify the classical crossings needed to unknot the diagram. Next to each such crossing, perform a virtual Reidemeister II move. Realize the new virtual crossings as shown in figure \ref{fig:unknotreal} to obtain a trivial realization.\qed

In the next section, we review representations of virtual links. We describe representations and their correspondence to virtual link diagrams and realizations. In particular, we examine the effect of homeomorphisms on a realization obtained from a representation.
In Section \ref{proofs}, we construct non-trivial realizations by applying a sequence of homeomorphisms to certain minimal genus representations.

\section{Representations and Realizations}

We recall
 \emph{representations} of virtual links \cite{kamada}, \cite{dk1} \cite{kup}. A representation of a virtual link $ \hat{L} $ is a pair $ (F,L)$, consisting of a link diagram on a closed, oriented two dimensional surface $ F $. Equivalence classes of representations are determined by orientation preserving homeomorphisms of the surface, handle addition and cancellation, and Reidemeister moves in the surface. Two representations related by a sequence of these moves are said to be \emph{stably equivalent}.

Orientation preserving homeomorphisms of the surface include deformations and Dehn twists \cite{purp}. To perform a Dehn twist on the surface, select a non-separating, simple closed curve on the surface. Cut the surface along this curve, to obtain a surface with two boundary components. Twist one of the boundary components one full rotation and reglue the boundary. This move does not induce knotting of the surface or change the genus.

Handle cancellation occurs along a \emph{cancellation curve}, an immersed copy of $ S^1 $ that does not intersect the immersed link. To cancel a handle, we cut the surface along the cancellation curve. This produces a surface with two boundary components. We then glue a disk to each boundary component. This procedure may separate the surface into two components. 

To add a handle to the surface, we identify a two disks on the surface that do not intersect the link diagram. We cut the surface along the boundary of these disks and remove the disks to form the surface $ F' $. Then we glue in a cylinder, $ S^1 \times I $, to the boundary of  $ F' $. The handle can be glued so that the surface becomes knotted.
It may be necessary to perform a sequence of Reidemeister moves in the surface prior to adding or removing a handle.

\begin{thm}
Equivalence classes of representations are in one to one correspondence with equivalence classes of virtual  link diagrams.
\end{thm}
\textbf{Proof:} See \cite{kamada}, \cite{kvirt}  \qed

\begin{rem}We can avoid the difficulty of performing Reidemeister moves before removing a handle. View (F,L) as an embedding of a link $ L $ into $ F \times I $. (See \cite{kup}). In this case, handle cancellation is performed along essential annuli in the surface.
Kuperberg used this approach in \cite{kup}.
\end{rem}

\begin{thm}[Kuperberg]Every stable equivalence class of links in a thickened surface has a unique irreducible representative.
\end{thm}
\textbf{Proof: } \cite{kup} \qed

This theorem proves that there is a unique minimal genus surface $ F $ among stably equivalent representations of a virtual link.

We examine the relationship between virtual link diagrams, representations, and realizations. We obtain
a representation of a virtual link diagram by the following process. Regard the virtual link diagram as a decorated immersion of $n$ copies of $ S^1 $ into the $ S^2 $. At each virtual crossing, select one arc in the crossing. Remove a small segment of this arc and attach a handle with an appropriately embedded arc to 
$ S^2 $.

 We can construct a representation $ (F,L)$ of a virtual link diagram that produces a specific realization when we \emph{forget} the underlying surface $ F $. 
For a fixed virtual link diagram with $ m $ virtual crossings, choose the over crossing arc for each virtual crossing. We then construct a representation $ ( F_m, L) $ (where $ F_m $ is an oriented, unknotted two dimensional surface of genus $ m $) by embedding the over crossing arc in the attached handle.  This process results in the selected realization. 

This representation can be reduced to a minimal genus representation via a sequence of homeomorphisms, Reidemeister moves in the surface and handle cancellations. However, this reduction may result in a trivial realization.

Given a representation of a virtual link that has minimal genus and an unknotted surface (and its corresponding realization), we 
can obtain a new realization of the link via a sequence of stably equivalent moves. Producing an equivalent  minimal genus representation with an unknotted surface restricts us to performing Dehn twists on the surface and Reidemeister moves in the surface. These moves do not change the genus or knot the surface $ F $.

\begin{figure}[htb] \epsfysize = 1.5 in
\centerline{\epsffile{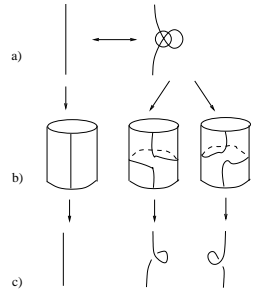}}
\caption{Virtual Reidemeister I Move}
\label{fig:vrIcorr}
\end{figure}

We examine how Dehn twists affect the virtual link diagram, representation, and realization.
In figure \ref{fig:vrIcorr}, we illustrate the virtual Reidemeister I move in row a). In row b) of this figure, we show representations of corresponding $ 1-1 $ tangles. The two representations are related by a Dehn twist applied to the meridian. In the last row, we remove the underlying surface and show the corresponding realizations. 

\begin{figure}[htb] \epsfysize = 1.5 in
\centerline{\epsffile{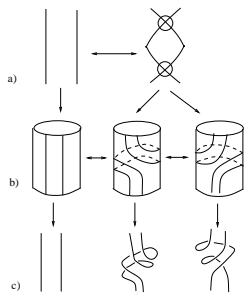}}
\caption{Virtual Reidemeister II Move}
\label{fig:vrIIcorr}
\end{figure}

We illustrate the effect of the virtual Reidemeister II move on a $ 2-2 $ tangle in row a) of figure \ref{fig:vrIIcorr}. In row b), we show two possible representations of this $ 2-2 $ tangle in a handle. (The representations are related by Dehn twists.) Finally, in row c), we again remove the underlying surface to obtain the realizations. In this case, the different realizations are not classically equivalent. 

For the remainder of this paper, twists (classical Reidemeister I moves) obtained by performing homeomorphism of the surface are suppressed in sketches of realizations. The twists do not affect the knotting or linking of the realization.

In the next section, we study virtual link diagrams with trivial realizations. Kishino's knot is an example of this type of virtual link. Every realization obtained from the standard diagram of Kishino's knot by choosing over/under markings is equivalent to the unknot.

\section{Constructing Non-Trivial Realizations} \label{proofs}
Let $ (F,L) $ be a representation of a virtual link $ \hat{L} $. The surface $F $ is the sum of $ n $ tori. Each handle in $ F $ has a merdian curve $ m $ and a longitude curve $ l $. 

In a minimal genus representation, some component $ L_i $ of $ L $ must intersect the longitude and meridan of each handle. If no component intersects a meridan (or longitude) then this curve is isotopic to a cancellation curve, contradicting our assumption that the genus was minimal. 

\begin{figure}[htb] \epsfysize = 1 in
\centerline{\epsffile{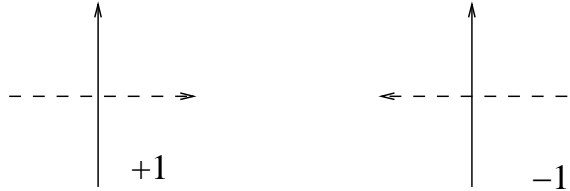}}
\caption{Computing the Intersection Pair}
\label{fig:intpair}
\end{figure}

We can measure the complicated relationship between the components of the link diagram $ L $  and the surface $ F $.
For each handle, choose an oriented meridian curve and an oriented longitude curve. We denote the meridian curve and the longitude of the $j^{th} $ handle as $ m_j $ and $l_j $ respectively. Assign an orientation to each component $ L_i $ of the link. For each handle and link component, we can compute the oriented intersection number between the component and the meridian and the compoment and the longitude.  
In figure \ref{fig:intpair}, the dashed line represents the meridian or longitude curve and the solid line represents the component. 
Each point of intersection contributes values as shown to the oriented intersection number.
A component that can be moved out of a handle through a sequence of Reidemeister II moves in the surface has an intersection number of zero.

We define the \emph{intersection pair} of a component and a handle. The intersection pair is denoted $ \lbrace p, n \rbrace $ where $ p $ is the oriented intersection number between the component and the meridian and $ n $ is the oriented intersection number between the component and the longitude.

\begin{thm} \label{int} Let $ \hat{L} $ be a non-trivial virtual link diagram. Let $ (F, L) $ be a minimal genus representation of $ \hat{L}$. If there is a component $ L_i $ with non-zero intersection pair, $ \lbrace p, n \rbrace $, then there exists a sequence of Dehn twists resulting in a stably  equivalent minimal genus representation with an unknotted surface corresponding to a non-trivial realization.
\end{thm}
 \textbf{Proof:}
Let $ (F, L) $ be a minimal genus representation of a virtual link $ \hat{L} $ with $ F $ unknotted. We assume that the realization is trivial (unknotted and unlinked).  We do not need to perform any Dehn twists if the realization is knotted or linked. By hypothesis, the component $ L_i $ and a handle in the surface have intersection pair $ \lbrace p , n \rbrace $. Without loss of generality, assume that $ p \neq 0 $.

 Take a small neighborhood ($ S^1 \times I $)  of the meridian curve. Isotope the link in the surface so that all oriented arcs (the orientation is inhereited from the link) in this neighborhood are of the form $ \star \times I $. Similarly, take a small neighborhood of the longitude and straighten the arcs in the neighborhood of the longitude. The dominant orientation of these two sets of arcs in is determined by the orientation of the majority of the arcs.
 
\begin{figure}[htb] \epsfysize = 1.5 in
\centerline{\epsffile{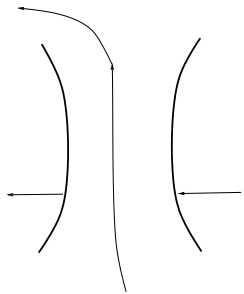}}
\caption{Original Representation}
\label{fig:1comp1}
\end{figure}

In the following figures,  we denote these collections of
arc by a single arrow pointing in the direction of the dominant orientation. We suppress the appearence of the other components. The sequence of Dehn twists may knot or link the other components of the diagram, but we have insufficient information to insure this conclusion. Our primary concern is to demonstrate that our process knots the component $ L_i $ in $ \mathbb{ R}^3 $.
In figure \ref{fig:1comp1}, we indicate schematically the position of $ L_i $.

\begin{figure}[htb] \epsfysize = 1.5 in
\centerline{\epsffile{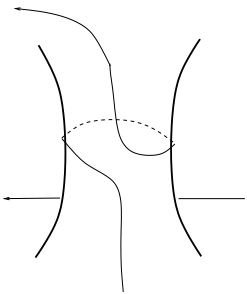}}
\caption{After Dehn Twist Along Meridian}
\label{fig:1comp2}
\end{figure}

Perform a Dehn twist along the meridian in the same direction as the oriented arc that intersects the longitude. As a result, we obtain the representation shown in figure \ref{fig:1comp2}.

\begin{figure}[htb] \epsfysize = 1.5 in
\centerline{\epsffile{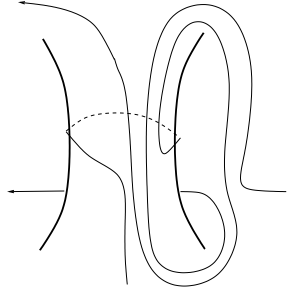}}
\caption{After Dehn Twist Along the Longitude}
\label{fig:1comp3}
\end{figure}

Dehn twist the representation along the longitude in the direction of the original arc that intersected the meridian. The result of this homeomorphism is pictured in figure \ref{fig:1comp3}.

\begin{figure}[htb] \epsfysize = 1.5 in
\centerline{\epsffile{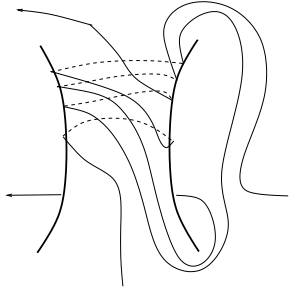}}
\caption{After Second Dehn Twist Along the Meridian}
\label{fig:1comp4}
\end{figure}

We perform a second Dehn twist along the meridian in the same direction (as the first Dehn twist). The result of this twisting is shown in figure \ref{fig:1comp4}

\begin{figure}[htb] \epsfysize = 1.25 in
\centerline{\epsffile{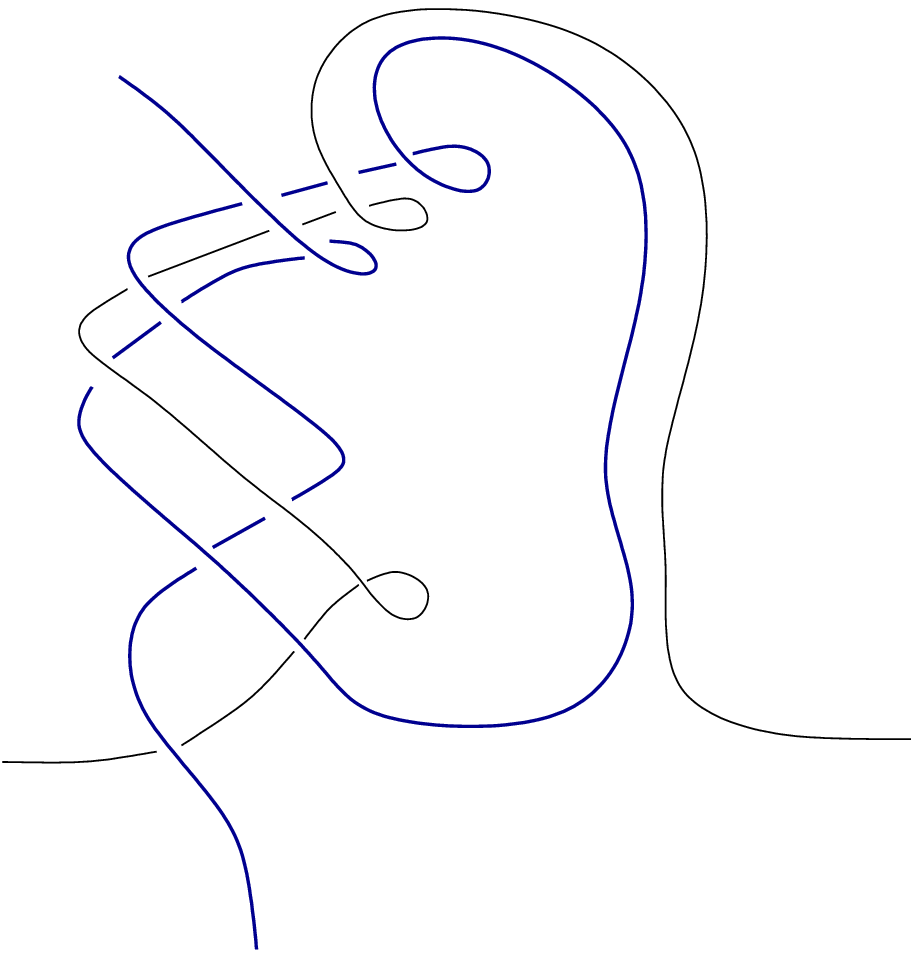}}
\caption{Realization obtained from figure \ref{fig:1comp4}}
\label{fig:twistsreal}
\end{figure}

We show a schematic of the realized tangle in figure \ref{fig:twistsreal}. The tangle contains a cabled trefoil.

\begin{figure}[htb] \epsfysize = 1.25 in
\centerline{\epsffile{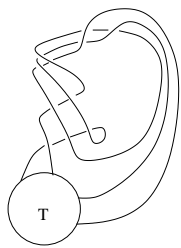}}
\caption{Schematic of the Realization}
\label{fig:realknot}
\end{figure}

A schematic of the realization is shown in figure  \ref{fig:realknot}.

If we can unknot the cabled trefoil in this schematic, each arc in the cabled trefoil must have a corresponding arc with opposite orientation. Then this pair of strands could be removed via a sequence of Reidemeister II moves. However, we assumed that there were $ p+k $ arcs with one orientation and $ k $ strands with the opposite orientation passing through the meridian. 
 \qed

\begin{rem}Note that if two or more components in the handle that have non-zero oriented intersection number with the meridian, a single Dehn twist along the meridian is sufficient to link the realization. 
Let $ L_j $ and $ L_i $ be a two components that have non-zero oriented intersection numbers $ p_j $ and $ p_i $ with a meridian curve. Since the realization obtained from $ (F,L) $ is trivial:
$ lk(L_i,L_j) = 0 $.
Perform a Dehn twist along the meridian and form the new minimal genus representation $ (F', L' ) $ with components $L'_i $ and $ L'_j $. 
 In the realization obtained from $ (F', L') $:
\begin{equation*}
lk (L'_i, L'_j)= \pm 2 p_i p_j.  
\end{equation*}
As result, the new representation results in a linked realization.
\end{rem}

If $ (F,L) $ is a representation of a virtual  and a single component intersects the meridian and the longitude of a handle exactly once, we may apply the following theorem. 

\begin{figure}[htb] \epsfysize = 0.75 in
\centerline{\epsffile{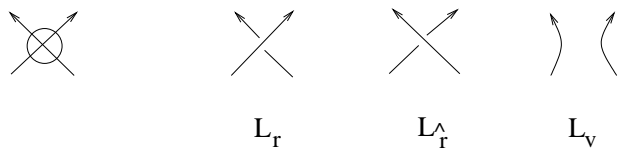}}
\caption{$ L_r$, $L_{\hat{r}} $, and $ L_v $}
\label{fig:choices}
\end{figure}

\begin{thm} \label{oneone}
Let $ L $ be a fixed virtual link diagram. If all realizations of this fixed diagram are trivial, we can construct a non-trivial realization by applying a single virtual Reideimeister II move.
\end{thm}
 \textbf{Proof:}
Let $ L $ be a fixed virtual link diagram such that all realizations are equivalent to the unknot. Let $ L_r $ be a realization obtained from $ L$. Let $ L_{\hat{r}} $ be a realization that differs from $ L_r $ at exactly one crossing, $ c$ . Without loss of generality, let $ sgn(c) =+1$ in $ L_r$. We illustrate $ L_r $, $ L_{\hat{r}}$, and $L_v$ in figure \ref{fig:choices}.

\begin{figure}[htb] \epsfysize = 1 in
\centerline{\epsffile{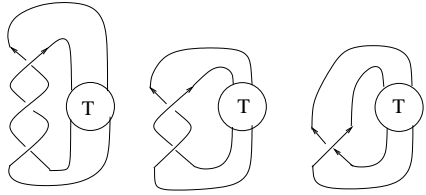}}
\caption{$ L_3 $, $ L_2 $, and $ L_r$}
\label{fig:triple}
\end{figure}

 Now,
\begin{equation*}
V(L_r)  = V(L_{\hat{r}} )= 1 
\end{equation*}
Note that if $ w(L_r)= -p+1$ then $ w(L_{ \hat{r} } ) =-p-1 $.
After orienting the diagram,
\begin{equation} \label{jskein}
A^4 V(L_r) - A^{-4} V(L_{\hat{r}}) = (A^{-2} -A^{2}) V( L_v)
\end{equation}
Hence, 
\begin{equation*}
\frac{A^4 -A^{-4}}{A^{-2} -A^{2}} =-( A^2 + A^{-2}) = V(L_v)
\end{equation*} 
Now, insert a virtual Reidemeister II twist in to the virtual link diagram $ L$. We construct a representation corresponding a realization of this diagram by performing a Dehn twist on the meridian of a handle. We obtain the realization $ L_3 $ from this representation. 
We now obtain the diagram $ L_3 $ as shown in figure \ref{fig:triple}.
We will also consider the diagram $L_2 $.
We now compute $ V(L_3) $. 
\begin{align*} 
A^4 V(L_3) &= (A^2 -A^{-2})V(L_2) + A^{-4} V(L_r) \\
V(L_3) &= (A^{-4})(A^{-2} -A^2) V(L_2) + A^{-8} V (L_r)
\end{align*}
We note that:
\begin{align*}
V(L_2) &= (A^{-4})(A^{-2} - A^{2}) V(L_r) + A^{-8} V(L_v) \\
V(L_2) &= (A^{-4})(A^{-2} -A^{2}) - A^{-8}(A^2 + A^{-2})
\end{align*}
Combining these terms:
\begin{equation*}
V(L_3)= -A^{-16} + A^{-12} + A^{-4}
\end{equation*}
As a result, $ L_3 $ is not equivalent to the unknot. \qed

We construct several examples using these results in the next section.

\section{Examples}

We present three examples utilizing the techniques discussed earlier in this paper.

\subsection{Kishino's knot}

\begin{figure}[htb] \epsfysize = 1 in
\centerline{\epsffile{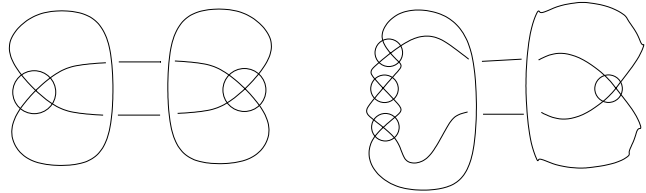}}
\caption{Kishino's Knot Example}
\label{fig:kishexa1}
\end{figure}

The standard diagram of Kishino's knot and an equivalent diagram is shown in figure \ref{fig:kishexa1}. Every realization obtained from the diagram on the left is classically equivalent to the unknot.

\begin{figure}[htb] \epsfysize = 1 in
\centerline{\epsffile{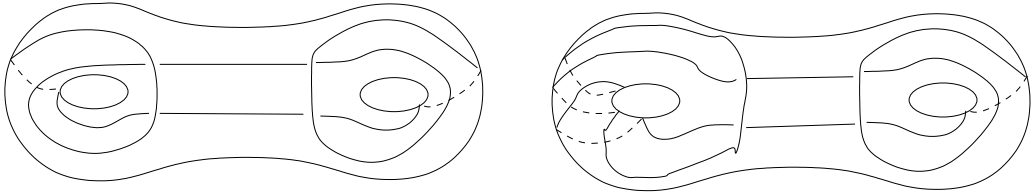}}
\caption{Two Representations of Kishino's Knot }
\label{fig:kishexa2}
\end{figure}

\begin{figure}[htb] \epsfysize = 1 in
\centerline{\epsffile{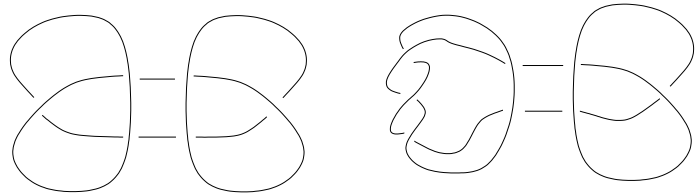}}
\caption{Realizations of Kishino's Knot}
\label{fig:kishexa3}
\end{figure}

In figure \ref{fig:kishexa2}, we show two equivalent representations of Kishino's knot. The two representations related by a Dehn twist about the meridian. 
The realizations obtained from these representations are illustrated in figure \ref{fig:kishexa3}. The classical knot diagram on the right is equivalent to the trefoil. 

\subsection{ Kauffman's Example }

\begin{figure}[htb] \epsfysize = 1 in
\centerline{\epsffile{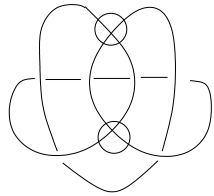}}
\caption{A Virtual Knot Diagram with Trivial Realizations }
\label{fig:kauffmanexa}
\end{figure}

The virtual knot diagram in figure \ref{fig:kauffmanexa} was suggested by Louis H. Kauffman. All realizations of this diagram are equivalent to the unknot.

\begin{figure}[htb] \epsfysize = 1.5 in
\centerline{\epsffile{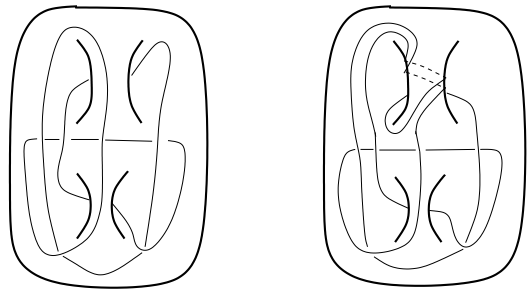}}
\caption{Two Equivalent Representations }
\label{fig:kauffmanexa2}
\end{figure}

\begin{figure}[htb] \epsfysize = 1.25 in
\centerline{\epsffile{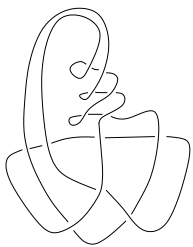}}
\caption{A knotted Realization}
\label{fig:kauffmanexa3}
\end{figure}

Two representations of this virtual knot diagram are illustrated in figure \ref{fig:kauffmanexa2}. Performing a Dehn twist along the meridian transforms the representation on the left into the representation on the right. 
This modified representation results the non-trivial realization shown in figure \ref{fig:kauffmanexa3}.

\subsection{An Example with Intersection Pair $ \lbrace 2,1 \rbrace $ }

\begin{figure}[htb] \epsfysize = 1 in
\centerline{\epsffile{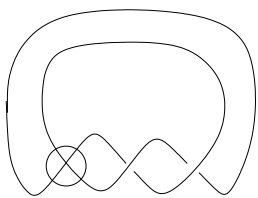}}
\caption{Virtual Knot Example}
\label{fig:2-1exa}
\end{figure}

The virtual knot diagram $ K $, shown in figure \ref{fig:2-1exa}, is non-classical. 

\begin{figure}[htb] \epsfysize = 1 in
\centerline{\epsffile{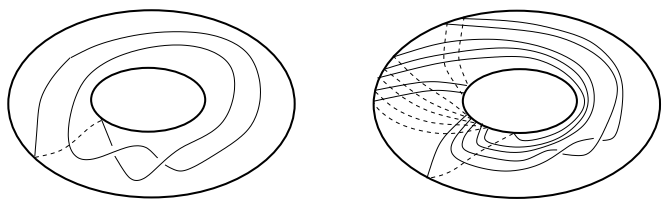}}
\caption{Example with Intersection Pair $ \lbrace 2, 1 \rbrace $}
\label{fig:2-1exa2}
\end{figure}

An unknotted, minimal genus representation is shown on the left in figure \ref{fig:2-1exa2}. This representation
results in a trivial realization and satisfies the hypothesis of Theorem \ref{int}. The intersection pair of this knot and torus is $ \lbrace 2, 1 \rbrace $. 

We apply the sequence of Dehn twists given in Theorem \ref{int} to obtain the equivalent representation pictured on the right of figure \ref{fig:2-1exa2}. 

\begin{figure}[htb] \epsfysize = 1.25 in
\centerline{\epsffile{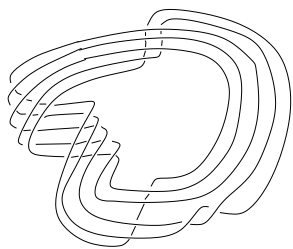}}
\caption{A Knotted Realization of $ K $.}
\label{fig:2-1exa3}
\end{figure}

This realization, $ K_r $, that is obtained from this representation is shown in figure \ref{fig:2-1exa3}.
The bracket polynomial of this knot is:
\begin{equation*}
\langle K_r \rangle = A^{-18} - A^{-14} - A^{-6} - A^2 + A^{22} + A^{30} + A^{38}
\end{equation*}
 
In conclusion, we make the following conjecture: 
\begin{conj}Let $ (F,L) $ be a minimal genus representation of a virtual link diagram. If every intersection pair is of the form $ \lbrace 0 , 0 \rbrace $ then there exists a minimal genus representation resulting in non-trivial representation.
\end{conj}

\end{document}